\documentclass[12pt]{article}
\usepackage{amsmath,amssymb}
\begin{document}

\section*{\centering Inhaltsverzeichnifs}

\section*{des zwei und zwanzigsten Bandes, nach den \centering Gegenst\"{a}nden.}

\section*{\centering I. Reine Mathematik.}
\section*{\centering 1. Analysis.}
\begin{tiny}
{Nr. der  \\  Abhandlung.} 
\end{tiny}
\hfill
\begin{tiny}
{Heft. seite}
\end{tiny}

\contentsline{section}
{\numberline {1.}\textbf{U}eber die  Aufl\"{o}sung der transcendenten Gleichungen.  Von Herrn Dr.\textit{M.A.Stern} zu G\"{o}ttingen. (Eine von der K\"{o}niglich-D\"{a} nischen Gesell-schaft der Wissenschaften gekr\"{o} nte Preisschrift.)}   {I.1}

\contentsline{section}
{\numberline {2.}Ueber die Zerlegung gebrochener algebraischer rationaler Functionen in Partialbr\"{u}che. Von Hrn. Prof.\textit{Oettinger} zu Freiburg i. Br.}{I. 63}

\contentsline{section}
{\numberline {5.}Schlufs dieser Abhandlung }{II 148}

\contentsline{section}
{\numberline {6.}Ueber die Bestimmung des Grades einer durch Elimination hervorgehen  den Gleichung. Von Hrn. Dr. \textit{Ferd. Minding} zu Berlin.}{II.178}

\contentsline{section}
{\numberline {7.}Ueber Transformation vielfacher integrale. Von Hrn. Dr.  \textit{Haedenkamp} zu Hamm in Westphalen.}{II. 184} 

\contentsline{section}
{\numberline{8.}Algemeine Aufl\"{o}sung der numerischen Gleichungen. Von Herrn Dr.Encke, Professor, Director der Sternwarte, Secretair der Akademie der Wissenschaften etc. zu Berlin.}{III. 193}

\contentsline{section}
{\numberline{9.}Einige Bemerkungen \"{u}ber die Mittel zur Sch\"{a}tzung der Convergenz  der allgemeinen Entwicklungs - Reihen mit Differenzen und Differentialen. Vom Herausgeber}{III.   249}
\begin{tiny}
{Nr.  der \\   Abhandlung.}
\end{tiny}
\hfill
\begin{tiny}
 Heft.   Seite.
\end{tiny}

\contentsline{section}
{\numberline {11.}De formatione et  proprietatibus Determinantium.
 Auct.  \textit{C. G. J. Jacobi,} prof. ord. math. Regiom } 
 { IV.   285}

\contentsline{section}
{\numberline {12.}De Determinantibus functionalibus. Auct.\textit{ C. G. J. Jacobi,}  prof. ord. math. Regiom }{ IV.   319}

\contentsline{section}
{\numberline{13.}De functionibus alternantibus earumque divisione  per productum e differentiis elementorum conflatum. Auct.\textit{ C. G. J. Jacobi,}  prof. ord. math. Regiom}{IV.   360}

\contentsline{section}
{\numberline{14.}Zur combinatorischen Analysis. Von  Hrn.\textit{ C. G. J. Jacobi,} Prof. ord. der Math. au der  Universit\"{a}t zu K\"{o}nigsberg in Pr }{ IV.   372}

\contentsline{section}
{\numberline{15.}Untersuchungen  \"{u}ber die Theorie der complexen Zahlen.Von Herrn Prof. \textit{ G.Lejeune Dirichlet} zu Berlin.  (Auszug aus einer der hiesigen  Akademie  der Wissenschaften am 27. Mai d. J. vorgelescnen Abhandlung.}{IV.   375}

 \subsection*\textit {iv \qquad Inhaltsverzeichnifs des zwei und zwanzigsten Bandes.}
              
\section*{\centering\numberline{2.} Geometrie.}

\contentsline{section}
{\numberline {3.}Elementare Ableitung eines zuerst von \textit{Legendre} aufgestellten Lehrsatzes der sph\"{a}rischen Trigonometrie. Von dem Herrn Hofrath, Prof. etc.Dr. \textsf{Gaufs} in G\"{o}ttingen}  {I.  96}

\contentsline{section}
{\numberline {10.}Geometrische Eigenschaftcn einer Factorentafel. Von Herrn Professor \textit{A. F. M\"{o}bius} in Leipzig.}{III.  276}

\contentsline{section}
{\numberline {16.}Eine Eigenschaft des Vierecks. Von Herrn Rechnungsrath \textit{Brune} zu Berlin.}{IV.  379}

 \section*{\centering \numberline{II.} Angewandte Mathematik}

    \section*{\centering  {Chronologie.}}

\contentsline{section}
{\numberline {4.}Zur Kirchenrechnung, Formeln und Tafeln. Von Hrn. Lic. \textit {Ferdinand Piper} in Berlin}      {ii.  97}

\contentsline{section}
 \numberline{} Druckfehler in diesem Bande

\newpage

\pagenumbering{arabic}

\section* {\centering 1.}

\section*{\centering Ueber die Aufl\"{o}sung der transcendenten}
   
\section* {\centering  Gleichungen.}

\section* {\centering \small\texttt (Vom Herrn Dr. \textit{M. A. Stern} zu G\"{o}ttingen.)}

\section*{\centering\small\texttt (Eine von der K\"{o}niglich-D\"{a}nischen Gesellschaft der  Wissenschaften gekr\"{o}nte Preisschrift.)}

\section*{\centering Vorwort.}

 \textbf{D}ie K\"{o}nigliche Gesellschaft der Wissenschaften zu Copenhagen hatte im Jahre 1837 folgende Preisfrage gestellt. 

Proponitur quaestio de aequationum transcendentium radicibus indagandiset quidem postulatur:

\begin{enumerate}

\item  	\textbf{U}t plene et perfecte deducantur interque se comparentur methodi ipsarum radices inveniendi, ita ut quaenam cuiuscunque sint virtutes quaenam imperfectiones accurate indicetur, quibusve casibus unaquaeque sit magis minusve accommodata.
    
\item       \textbf{U}t diligenter inquiratur quatenus vel quibus saltem adhibitis cautionibus methodos, quibus vulgo in algebraicis aequationibus radices reales aut ab imaginariis separentur aut inter se, ad transcendentes quoque extendere liceat.

\item  	\textbf{U}t exponatur conspectus,quantum fieri possit, plenus tam specialium aequationum quam generum earum, quae quidem forma transcendenti in gravissimis analyseos applicatae partibus occurrunt, simul cum regulis, quin fortasse tabulis ad usum ipsum accomodatis, quibus revera faciliores ac breviores reddantur calculi illi radicum, alias saepe prolixissimi.

\end{enumerate}

\footnote {Crelle's Journal  d. M. Bd. XXII. Hft I.}      

  	\textbf{D}as Folgende ist ein genauer Abdruck der Schrift, die ich der   K\"{o}nigl. G. d. W. im December 1837 \"{u}berreicht habe; ich habe mir nur einige aufserwesentliche Aenderungen erlaubt, wie namentlich, dafs ich da,  wo ich, der Bestimmung der Schrift gem\"{a}fs, von meinen eigenen Arbeiten  als denen eines Dritten  sprechen mufste, dies nun ge\"{a}ndert habe.\\

\hspace{2em}\textbf{G}\"{o}ttingen, den 5. Januar 1840.

\section*{\centering 6.}                            
\section* {\centering  {Ueber die Bestimmung des Grades einer durch Elimination hervorgehenden Gleichung.}} 

 \section* {\centering \small\texttt (Von Herrn Dr. \textit{Ferd. Minding} zu Berlin.)}

\textbf{E}s sind zwar verschiedene Methoden bekannt, um aus zwei algebraischen Gleichungen, zwischen zwei Unbekannten, eine neue, nur noch eine Unbekannte enthaltende Gleichung herzuleiten; h\"{a}ufig aber w\"{u}nscht man nur den Grad dieser Endgleichung zu wissen, nicht sie selbst aufzustellen, und wenn ich nicht irre, so ist eine Regel, welche diesen Grad mit der einem so elementaren Gegenstande angemessenen Leichtigkeit finden lehrte, bis jetzt nicht gegeben worden. Man weifs zwar, wenn die Gleichungen beziehungsweise von den Graden $h$  und $k$  sind, d. h. wenn die h\"{o}chste Summe der Exponcnten eines Gliedes in der einen $h$, in der anderen $k$ betr\"{a}gt, dafs alsdann der gesuchte Grad dem Producte $hk$ h\"{o}chstens gleichkommen kann; hiermit ist aber nur eine Grenze gegeben, von welcher der wirkliche Grad oft sehr abweicht. Um diesen zu finden, ist vor allemn\"{o}thig, die wahre Form der Endgleichung festzustellen. Man habe folgende Gleichungen:

\begin{alignat*}{2}
1.\qquad{f}(x,y)&=A_0y^m + A_1y^{m-1} + A_2y^{m-2}+\dots+A_{m-1}y+A_m &= 0, \\
2.\qquad \theta(x,y)&=B_0y^n+B_1y^{n-1}+B_
2y^{n-2}+\dots+B_{n-1}y+B_{n} &= 0,
\end{alignat*}
in welchen  die Buchstaben $A$ und $B$ mit angeh\"{a}ngten Zeigern beliebige ganze Polynome in $x$ bedeuten. L\"{o}set man die Gleichung (2.) nach  $y$  auf, bezeichnct ihre Wurzeln mit   ${y_1, y_2,\dots y_n,}$ und bildet das Product

\begin {equation*}                        
3.\qquad P = f (x, y_1) \cdot f (x, y_2)\dots f(x,y_n),
\end{equation*}
so ist  ${B_0^m}\cdot{P}$ eine ganze Function von $x$, und wenn man setzt 

\begin{align*}
4.\qquad {B_0^m}\cdot{P}&=\psi{x},  \\
{\text{so  ist}}\qquad\qquad    5.\qquad \psi{x}&=0
\end{align*}
die verlangte Endgleichung.

     \textbf{U}m zu zeigen, dafs $B_0^m\cdot{P}$ eine ganze Function ist, bemerke man zuerst, dafs $P$ eine rationale Function von $x$ ist, welche, wenn die darin vorkommenden symmetrischen Functionen von ${y_1, y_2,\dots y_n,} $  vermittelst der Gleichung $(2.)$ in $x$ ausgedr\"{u}ckt  werden, nur eine Potenz von $B_o$ zum Nenner erhalten kann.  Bezeichnet man irgend eine  jener symmetrischen Functionen durch
\begin{equation*}
S = {y_1^{m_1}} {y_2^{m_2}}{y_3^{m_3}}\dots{y_n^{m_n}}+\dots,
\end{equation*}
wo die nachfolgenden durch Puncte angedeuteten Glieder aus dem ersten durch Verwechselung von  $y_1, y_2,\dots y_n$  entstehen, so kann keiner der  Exponenten $m_1, m_2, \dots m_n$  gr\"{o}fser sein als $m$.Man setze nun:

\begin{align*}
{S_1}&=\frac{1}{y_1^{m-{m_1}}\cdot{y}_2^{m-{m_2}}\dots{y}_n^{m-{m_n}}}+\dots, \\
\hspace{-15em}{\text{mithin}}\qquad\quad\qquad\;\;{S}& = {(y_1,y_2\dots y_n)}^m {S_1}=(-1)^{mn}\frac{B_n^m}{B_o^m}{S_1}.
\end{align*}
Das Zeichen $S_1$ bedeutet eine ganze symmetrische Function der ungekehrten Wurzeln von $(2.)$ mithin ist der Werth von $S_1$ 
ein rationaler Bruch, der nur eine Potenz von $B_n$  zum Nenner haben kann; setzt man daher $S_1=\frac{Z}{B_{n}^\lambda},$
wo Z ein ganzes Polynom in $x$ oder genauer eine ganze Function der Polynome
${B_0, B_1, \dots {B}_n}$ und $\lambda$ eine positive ganze Zahl ist, so kommt

\begin{equation*}
S=(-1)^{mn}\frac{B_{n}^{m}\cdot{Z}}{B_{0}^{m}\cdot{B}_{m}^{\lambda}}
\end{equation*}
Da nun $S$ offenbar nur eine Potenz von $B_0$ zum Nenner haben kann, so mufs $B_{n}^{\lambda}$ in dem vorstehenden Z\"{a}hlcr aufgehen, mithin ist $B_{0}^{m}{S}$ und folglich auch $B_{0}^{m}{P}=\psi{x}$ eine ganze Function. 
                   
         Bezeichnet man die Wurzeln der Gleichung $(1.)$ nach  $y$ aufgel\"{o}st,  durch  $\eta_{1}, \eta_{2}, \cdots\eta_{m},$ setzt

\begin{equation*}
{Q}=\theta(x,\eta_{1})\cdot\theta(x,\eta_{2})\cdots\theta(x,\eta_{m})
\end{equation*}
und bemerkt dafs

\begin{equation*}
\theta(x,\eta_1)=B_0(\eta_1-y_1)(\eta_1-y_2)\dots(\eta_1-y_n),\text {u.s.f.},
\end{equation*}
so wird: 
\begin{multline*}
Q=B_0(\eta_1-y_1)\dots(\eta_1-y_n)\times{B}_0(\eta_2-y_1)\dots(\eta_2-y_n)\times\dots\\
\dots\times{B}_0(\eta_m-y_1)\dots(\eta_m-y_n),
\end{multline*}
und weil
\begin{equation*}
A_{0}(y_1-\eta_1)(y_1-\eta_2)\cdots(y_1-\eta_m)=f(x,y_1),\text{u.s.f.},
\end{equation*}
so folgt:

\begin{equation*}
6.\qquad     A_{0}^{n}Q=(-1)^{mn}B_{0}^{m}P.
\end{equation*} 

\textbf{H}ieraus ergiebt sich, dafs $\psi{x}$=0 die verlangte Endgleichung ist. N\"{a}mlich f\"{u}r jeden der Aufgabe zusagenden Werth von $x$ wird nothwendig  $P = 0$ (eben so auch $Q = 0$ ),
 also $\psi{x}$=0. Sollte ferner diese Gleichung einen \"{u}berfl\"{u}ssigen Factor enthalten, so w\"{a}re f\"{u}r einen solchen $\psi{x}$ = 0, zugleich aber weder $P$ = 0 noch $Q$ = 0; alsdann m\"{u}fsten, wegen $(4.)$ und $(6.)$, $A_0$ und $B_0 $ zugleich verschwinden, was im Allgemeinen  nicht m\"{o}glich ist. Haben in einem besonderen Falle  $A_0$ und $B_0 $ einen Factor gemein, so ist allemal  auch das Polynom $\psi{x}$ durch diesen Factor theilbar, weil cs immer, wie leicht zu sehen, von folgender Form ist: $\psi{x}=A_0U+{B_0}V$, 
in welcher $U$ und $V$ ganze Polynome sind; da man jedoch einen solchen Fall stets durch eine unendlich kleine Aenderung der Coefficienten beseitigen kann, und zwar ohne den Grad eines derselben zu  \"{a}ndern, so folgt, dafs die Gleichung $\psi{x} = 0$ in keinem Falle eine der Aufgabe fremde Wurzel darbietet. 
 
Der Grad des Polynoms $\psi{x}$ ergiebt sich nun auf folgende Weise. Man hat
\begin{equation*}                 
\psi{x}=B_{0}^{m}\cdot{f}(x,y_1)\cdot{f}(x,y_2)\dots{f}(x,y_n).
\end{equation*}
Entwickelt man die Wurzeln  ${y}_1,  {y}_2,\dots {y}_n$  der Gleichung $(2.)$ nach fallenden Potenzen von $x$, und setzt die erhaltenen Reihen anstatt jener in vorstchenden Ausdruck, so werden alle gebrochnen und negativen Potenzen von $x$ sich gegenseitig aufheben und das Polynom $\psi{x}$ wird unver\"{a}ndert, wie vorhin, hervorgehen. Da nur der Grad von $\psi{x}$ verlangt wird, so setze man statt jener Reihe nur ihre ersten Glieder, die f\"{u}r ${y}_1, {y}_2\dots{y}_n$  beziehungsweise sein m\"{o}gen: $c_1{x}^{h_1},
 {c}_2{x}^{h_2}\dots {c}_n{x}^{h_n}$.  Das Verfahren, durch welches die Reihen und namentlich die h\"{o}chsten Exponenten $h_1, {h}_2, \dots {h}_n$ oder die Grade der Wurzeln gefunden werden, ist hinl\"{a}nglich bekannt; man vergleiche z. B. \textit{Lacroix Traite} S. 223 der  ersten Ausgabe, wo die Entwickelung nach steigenden Potenzen gezeigt wird . Man bestimme hierauf den	h\"{o}chsten Exponenten von $x$ in jeder der Functionen  $f(x, {c}_1{x}^{h_1}), f(x, {c}_2{x}^{h_2}),\dots $   oder die \textit{Grade} der  Functionen  $f(x, y_1), f(x, y_2), \dotsc,$ 
 welche mit ${k}_1,  {k}_2,\dots {k}_n$  bezeichnet werden m\"{o}gen. Diese k\"{o}nnen ganz oder gebrochen, aber nie negativ sein, weil $A_n$  wenigstens vom Grade Null ist. Wird endlich noch der Grad von $B_0$  mit $b$  bezeichnet, so ist 
\begin{equation*}
7.\qquad  mb+{k}_1 +  {k}_2 +{k}_ 3 + \dots + {k}_n 
\end{equation*}
nothwendig eine \textit{ganze}  Zahl, welche den h\"{o}chsten 
Exponenten von $\psi{x}$ oder den gesuchten Grad der Endgleichung angiebt.
In besonderen F\"{a}llen kann man noch die Werthe von $c_1,  {c}_2,\dots {c}_n$  ber\"{u}cksichtigen, um zu sehen, ob der Coefficient des h\"{o}chsten Gliedes in einem der Factoren  $f(x,y_1), \dots$  von $\psi{x}$ und mithin in $\psi{x}$ selbst vielleicht gerade Null  wird, und in  einem solchen Falle wird man gen\"{o}thigt sein, auch die folgenden Glieder der Reihen f\"{u}r  ${y}_1, {y}_2, \dots {y}_n$  
 theilweise in Rechnung zu bringen; es wird  jedoch nicht erforderlich sein diese Andeutung hier weiter auszuf\"{u}hren; vielmehr ist klar, dafs im Allgemeinen der obige Werth $(7.)$ den wirklichen Grad des Polynoms $\psi{x}$  darstellt.
		         
\textbf{E}s seien z. B. folgende zwei Gleichungen gegeben, in welchen das Zeichen ${(x)}^{\mu}$ ein Polynom in $x$ vom Grade $\mu$ anzeigt:

\begin{align*}
f(x,y)& = (x^2)y^4+(x^2)y^3+(x^4)y^2+(x^5)y+(x^5) = 0,\\
\theta(x, y)& = (x^8)y^5+(x^6)y^4+(x^9)y^3+(x^4)y^2+(x^3)y+(x^4) = 0.
\end{align*}
Diese Gleichungen sind von \textsf{Gten} und vom 13ten Grade, also ist der Grad der Endgleichung nicht h\"{o}her  als $6.13 = 78$. Um ihn genau zu finden, berechne man die Grade der Wurzeln $y$ von $Q = (x, y) = 0$;  man findet sofort $h_1 = h_2 = \frac{1}{2}, h_3 = h_4 = h_5 = {-} \frac{5}{3}$.  Hieraus folgen die Grade von $f(x, y_1), \dots,$   n\"{a}mlich $k_1 = k_2 = \frac{11}{2}, k_3 = k_4 = k_5 = 5$;         ferner ist $B_0 = (x)^8$, also $b = 8$, und $m = 4$, also der Grad der Endgleichung $mb + k_1 + k_2 + k_ 3 + k_4 +k_5 = 4.8 + 11 + 15 =58$.

    Wenn man die gegebenen Gleichungen, anstatt nach $y$, nach $x$  ordnet, um nach obiger Regel  den Grad der Endgleichung in $y$ zu  suchen, so findet man nicht immer denselben Werth f\"{u}r  diesen wie f\"{u}r den vorigen Grad. Zur Erkl\"{a} rung dieses Umstandes mufs man bemerken, dafs die Endgleichung in $x$ nur die \textit{endlichen} Werthe von $x$ ergiebt, welche beiden vorgelegten Gleichungen zu gen\"{u}gen gceignet sind. Steigt also die Endgleichung in $y$ auf einen h\"{o}heren Grad als die in $x$, so geh\"{o}ren nothwendig einige  Wurzeln der Gleichung in $y$  zu uncndlichen Wcrthen von $x$. Es ist auch allemal leicht, diese Werthe durch eine unendlich kleine Aeuderung der Coefficienten in einer der vorliegenden Gleichungen  zum Vorschein zu bringen, und die Ungleichheit der Grade der Endgleichungen  zu tilgen. Es seien n\"{a}mlich die Gleichungen $(1.)$ und $(2.)$ nach $x$ geordnet, folgende:

\begin{eqnarray*}
f(x,y)& =& \alpha_0{x}^\mu+\alpha_1{x}^{\mu-1}+ \dots+\alpha_\mu = o,\\     
\theta(x, y)& =& \beta_0{x}^\nu+\beta_1{x}^{\nu-1}+ \dots+\beta_\nu = o,
\end{eqnarray*}
wo $\alpha_0, \alpha_1, \dots, \beta_0, \dots\beta_\nu$  ganze Polynome in $y$ sind. Wenn nun weder  $A_0$ mit  $B_0$, noch $\alpha_0$, mit $\beta_0$  einen Factor gemein hat, so k\"{o}nnen weder unendliche Werthe von $y$ f\"{u}r endliche von $x$, noch unendliche von $x$ f\"{u}r endliche von $y$ Statt finden; folglich kann alsdann zwischen den Graden der Endgleichungen in $x$ und in $y$  kein Unterschied sein.  Man braucht also,wenn gemeinsame Factoren  zwischen $A_0$  und $B_0$  oder  $\alpha_0$  und  $\beta_0$ vorhanden  sind, nur einen Coefficienten in $A_0$  und einen  in  $\alpha_0$  zu \"{a}ndern, um f\"{u}r die Endgleichungen in $x$ und $y$ gleiche Grade zu erhalten. Setzt man nachher diese  Aenderungen  gleich Null, so kann man die Coefficienten  der h\"{o}chsten Glider  der Endgleichungen pr\"{u}fen, um zu entscheiden,  wie  viele Werthe von $x$   und wie viele  von $y$  unendlich werden, und wie viele \textit{endliche} L\"{o}sungen der vorgelegten Gleichungen schliefslich vorhanden sind. Diese Ausf\"{u}hrung der Rechnung ist jedoch unn\"{o}thig, wenn man die vorgetragene Regel  geh\"{o}rig anwendet.  Man habe z. B. folgende Gleichungen:

\begin{eqnarray*}
f(x, y) & = &(a+bx^2)y^4+(c+ex)y^2+gx^3y+h +kx^2+lx^3 = 0, \\
\theta(x,y) & = & \beta{x}^5y^2+(\gamma+\delta{x}^2)y+ \lambda+\mu{x}^4 = 0, 
\end{eqnarray*}
oder nach $x$ geordnet:

\begin{eqnarray*}
f(x,y) & = &(l+gy)x^3 + (k+by^4)x^2+ey^2x +h + cy^2 +ay^4 = 0, \\
\theta(x, y) & = &\beta{y}^2x^5+ \mu{x}^4 +\delta{yx}^2 +\lambda + \gamma{y} = 0.
\end{eqnarray*}
Hier ist $A_0 = a+ bx^2, B_0 =\beta{x}^5, \alpha_0 = l+gy, \beta_0 = \beta{y}^2$;\quad folglich haben, wenn weder  $a = 0$, noch ${l} = 0$, $A_0$ und $B_0$, so wie $\alpha_0$ und $\beta_0$ keinen  gemeinsamen Factor, daher die Grade der Endgleichungen in $x$  und $y$  \"{u}bereinstimmend gefunden werden $= 26$. Setzt man aber zugleich $a = 0 $ und  $\textit{l} = 0 $,  und berechnet alsdann die Grade der   Endgleichungen, so findet man $25$ f\"{u}r die Gleichung in $x$, und $24$ f\"{u}r die in $y$. Durch das Verschwinden von $a$  und $l$ werden also zwei Wurzeln der vorigen  Endgleichung in $y$ und eine der vorigen Endgleichung in $x$  unendlich; zugleich aber wird auch die neue Endgleichung in $x$ durch ${x}^2$, den gemeinsamen Factor von $A_0$  und $B_0$, so wie die neue Endgleichung in $y$ durch $y$, den gemeinsamen Factor von $\alpha_0 $ und $\beta_0$,  theilbar. Von den $26$ endlichen L\"{o}sungen, wclche den anf\"{a}nglichen Gleichungen zukamen, bleiben also $23$ im Allgemeinen  noch endlich, wenn $a$  und $l$ verschwinden;  die drei \"{u}brigen hingegen sind: $x_{24} = 0$, $y_{24} = \infty$;  $x_{25} = 0$, $y_{25} = \infty$; $x_{26} = \infty$, $y_{26} = 0$. Man sieht, wie hier die gesuchte Anzahl der endlichen L\"{o}sungen durch wiederholte Anwendung der vorgetragenen Regel und Vergleichung der Resultate gefunden wird.

\textbf{N}achschrift. Nach Beendigung des Vorstehenden ist mir ein neues Werk zu Gesicht gekommen: ;,,\textbf{S}ystem der Algebra von Dr. \textit{P.J.E. Finck,} Professor zu Strafsburg; Leipzig bei Barth, 1841,'` in welchem S. 405 zur Bestimmung des Grades der Endgleichung eine viel genauere Regel angegeben wird,  als diejenige, deren im Eingange des vorstehenden Aufsatzes erw\"{a}hnt ist. Die Regel ist, dem Inhalte nach, folgende: Wenn alle Coefficienten  $A_0, A_n, \dots A_m$ der Gleichung $(1.)$ vom Grade  $m^\prime$ und alle Coefficienten $B_0, B_1, \dots B_n$ der Gleichung $(2.)$ vom Grade $n^\prime$ sind,
so ist der Grad des Polynoms $\psi{x}$, welches die Endgleichung liefert, folgender: ${m}n^\prime + n m^\prime$. Der Beweis dieses  Satzes, der zwei Seiten des Lehrbuches f\"{u}llt, folgt sehr einfach aus dem Obigen; denn in diesem Falle sind alle $y$ aus der Gleichung $(2.)$ vom Grade 0, d.i. $h_1 = h_2 = \dots = h_n = 0$, folglich $k_1 = k_2 = \dots = k_n = m^\prime$; zugleich ist $b 
= n^\prime$, weil $B_0$ vom Grade $n^\prime$; folglich der Grad der Endgleichung: $m b + k_1 + k_2 + \dots + k_n  =m n^\prime + n m^\prime$, w.z.b.w
Wird diese Regel auf F\"{a}lle angewendet, in welchen die Coefficienten von ungleichen Graden sind, so ist das Resultat nicht mehr zuverl\"{a}ssig, weil es die Ausgleichung jener Grade voraussetzt, welche fremdartige Wurzeln herbeif\"{u}hrt. Bei der im gegenw\"{a}rtigen Aufsatz vorgetragenen Regel werden dagegen, um den Grad von $\psi{x}$ zu finden, die Grade der Coefficienten  nur so in Rechnung gebracht, wie sie gegeben sind.

\end{document}